\documentclass[a4paper,12pt,leqno]{article}
\usepackage[T1]{fontenc}
\usepackage[utf8]{inputenc}

\usepackage{mathptmx}       
\usepackage{helvet}         
\usepackage{courier}        
\usepackage[bottom]{footmisc}
\usepackage{url}
\usepackage{latexsym}
\usepackage{amssymb,amsmath}
\usepackage{amsfonts,amsthm}
\usepackage{color}
\usepackage[a4paper,margin=1in]{geometry}

\theoremstyle{plain}

\newtheorem*{theorem*}{Theorem}
\newtheorem*{lemma*}{Lemma}
\newtheorem*{proposition*}{Proposition}
\newtheorem*{corollary*}{Corollary}

\theoremstyle{definition}

\newtheorem*{remark*}{Remark}
\newtheorem*{definition*}{Definition}
\newtheorem*{example*}{Example}

\begin{document}
\title{Feller's Contributions to  Mathematical Biology}
\author{
    Ellen Baake\thanks{Faculty of Technology, Bielefeld University, Box 100131, 33501 Bielefeld, Germany} \, and Anton Wakolbinger\thanks{Institut f\"ur Mathematik, Goethe-Universit\"at, Box 111932, 60054 Frankfurt am Main, Germany}   
}
\date{}
\maketitle
\begin{abstract}
This is a review of William Feller's important contributions to mathematical biology. The seminal paper \cite{Feller 1951d} {\em Diffusion processes in genetics} was particularly influential on the development of stochastic processes at the interface to evolutionary biology, and interesting ideas in this direction (including a first characterization of what is nowadays known as ``Feller's branching diffusion'')  already shaped up in the paper \cite{Feller 1939a} (written in German) {\em The foundations of a probabistic treatment of Volterra's theory of the struggle for life}. Feller's article {\em On fitness and the cost of natural selection}  \cite{Feller 1967c} contains a critical analysis of the concept of {\em genetic load}.

The present article will appear in: {\em Schilling, R.L., Vondracek, Z., Woyczynski, W.A.: The Selected Papers of William Feller. Springer Verlag}.
\end{abstract}

\section{Introduction} \label{sec:BW-intro}
Feller had a persistent  interest in biology. This is documented in numerous
examples from mathematical genetics in his monograph
\cite{Feller 1950,Feller 1966d},
and by a couple of  influential research papers at the interface
of population biology and probability theory. Looking back at these
papers in historical perspective is highly rewarding:
They are cornerstones of biomathematics;
they mirror the development of probability theory of their time;
and at least one of them (\cite{Feller 1951d}) had lasting impact on probability theory.

Feller's  important papers on the interface to biology are \cite{Feller 1939a},
\cite{Feller 1951d}, and \cite{Feller 1967c}. The first one addresses
general population dynamics, the other two are mainly
concerned with  models in population genetics.
The area of \emph{population dynamics} is concerned with the growth,
stabilisation, decay, or extinction of populations. Models of population
dynamics describe how the \emph{size} of populations changes over time
under given assumptions on birth and death rates of individuals,
which may depend on the current population size since individuals
interact (e.\,g.~compete) with each other. In contrast,
\emph{population genetics} is concerned with the \emph{genetic composition}
of populations under the action of various evolutionary processes,
such as  mutation  and selection. Naturally, there is no sharp
boundary between the fields, as we will also see in Feller's contributions.
Let us now look at them.

\section{Feller and population dynamics}
In \cite{Feller 1939a}, a paper still in German entitled (in English
translation) \emph{The foundations of a probabilistic treatment
of  Volterra's theory of the struggle for life}, Feller presents
a synthesis of two fundamental developments that both
started in 1931. On the one hand, Volterra presented his book
\emph{Lessons about the mathematical theory of the struggle for life}
\cite{BW-Volterra1931};
on the other hand, Kolmogorov published his seminal
paper \emph{On analytical methods in probability theory} \cite{BW-Kolmogorov1931}.
Volterra's book laid the
foundations for the \emph{deterministic} description of population dynamics in
terms of systems of ordinary differential equations that model
birth, death, and interaction of individuals. These models imply
that populations are so large that random fluctuations can be
neglected, and population sizes are measured in units so large
that the size can be considered a continuous quantity.
Kolmogorov presented
the general and systematic
formalism for the description of  \emph{stochastic} dynamics
in terms of Markov chains in continuous time; in particular, he found the
description for the evolution of probability weights and the transport
of expectations in terms of differential equations,
which we know today under the names of \emph{Kolmogorov forward equations}
and \emph{Kolmogorov backward equations}.

In his 1939 paper, Feller ties  these two fundamental developments together
by applying Kolmogorov's new formalism to some examples of Volterra's
population dynamics. We see here the birth of the stochastic description
of population dynamics, which today has its firm place in mathematical
biology,
and is highly developed both in analytical terms and in terms of
simulations.

Feller's paper is devoted to the description of single populations
(except from a small excursion to predator--prey models in the end)
and consists of two large parts. The first establishes
the Kolmogorov forward equations (KFE) for the Markov jump
processes (namely, birth-and-death processes)
that describe finite populations (remarkably, there is no mention
of the Kolmogorov backward equations in this paper).
The second part discusses a continuum analogue of such processes, a special case of which seems to be the first appearance of what
today is called \emph{Feller's branching diffusion}.

It is remarkable to see (and a pleasure to read) that Feller notices
some of the crucial relationships between corresponding deterministic
and stochastic models in this early paper, which appear as a central
theme.

For the sake of clarity,
let us make explicit here the two fundamental limits of birth-death processes that are addressed in \cite{Feller 1939a}.
Consider a birth-death process $K_N(t)$ with  birth rate $n \lambda$ and
death rate $n \mu$ when in state $n$, with $N$ being the  initial population size.  Then, as the initial population size $N$ tends to $\infty$, the sequence of process $(K_N(t)/N)_{t\ge0}$, $N=1,2,\ldots$, converges in distribution to the solution of the differential equation
\begin{equation}\label{eq:BW-linode}
\dot x = (\lambda-\mu)x, \quad x(0) =1.
\end{equation}
This reflects a dynamical version of the Law of Large Numbers (see e.\,g.~\cite{BW-Kurtz1970}). (Notably, due to the linearity, the expectation $M(t):= \mathbb E(K_1(t))$ satisfies \eqref{eq:BW-linode} as well.) A different kind of limit emerges if one assumes that the individual split and death rates $\lambda$ and $\mu$ depend on $N$ and the process is \emph{nearly critical} in the sense that 
$$
    \lambda_N= \beta + \theta_1/N
     \quad\text{and}\quad
     \mu_N= \beta+ \theta_2/N, 
$$
with $\theta_1-\theta_2 =:\alpha$.  The Law of Large Numbers then says that the limit of the processes $(K_N(t)/N)_{t\ge0}$ is the constant $1$. However, on a larger time scale the fluctuations become visible: the sequence of processes $(K_N(Nt)/N)_{t\ge0}$ converges in distribution to the solution of the \emph{stochastic} differential equation \eqref{eq:BW-51} stated in paragraph \ref{BW-twoSDEs}, whose diffusion equation is \eqref{eq:BW-Feller_branching}. This is a prototype of a diffusion limit for birth-death processes.
In  \cite{Feller 1939a},  these limiting procedures are not made explicit (but see \cite{Feller 1951d} for a major step in this direction). Feller in 1939 goes rather the other way,
in search for stochastic processes that correspond to a given
deterministic model.
Let us now explain the major lines of his article.

\subsection{Markov jump processes for population dynamics}
In the first part (Sections 1--4),
devoted to the stochastic description of finite
populations, Feller explains a variety of birth-and-death processes
and sets up the Kolmogorov forward equations for them, i.\,e.~he establishes
 the system of differential equations that describe
how the probability weights for the number of individuals
alive at time $t$ evolve over time. He starts with the simple
linear death process (where each individual dies at rate $\lambda$,
independently of all others), proceeds via the corresponding birth
process and the linear birth-and-death process and finally arrives
at the general birth-and-death process. In an individual-based picture, the latter includes interaction
between individuals, so that the birth and/or death rates
are no longer linear in the number of individuals.
The case of logistic growth, which includes a quadratic competition
term, serves as an important example; the case of `positive interaction'
(such as symbiosis) is not treated explicitly here. Let us comment on the major
insights of this part.

\subsubsection{Kolmogorov equations, their solutions, and relationship
with deterministic description.}
Feller notices that for a given net reproduction rate $\alpha$
per individual, by choosing $\lambda -\mu = \alpha$, one obtains a variety of linear birth-death processes whose expectation value $M(t)$ satisfies one and the same  ODE \eqref{eq:BW-linode}, whereas for $\alpha > 0$, there is exactly one linear pure birth process ($\lambda = \alpha, \mu = 0$) with this property.
Feller states this
ambiguity explicitly when discussing  \emph{logistic growth}.
Its deterministic
version is given by the differential equation
\begin{equation}\label{eq:BW-logistic}
  \dot m = m (\lambda - \gamma m) =: f(m),
\end{equation}
which Feller also calls the \emph{Pearl-Verhulst} equation.
Here $m$ is shorthand for $m(t)$, the `deterministic version' of the population size at time $t$,
$\lambda$ denotes the  per capita
net reproduction rate in the absence of competition, and $\gamma$
is the competition parameter. Again, Feller
notices that there are many possibilities in terms of
birth-death processes that correspond to \eqref{eq:BW-logistic}. They are parametrised in his Eq.~(27),
which describes the
process with per capita birth at rate $\omega - \nu n$ and per capita death at rate
$\tau - \sigma n$ if there are currently $n$ individuals. 
Here, we have renamed $\gamma$ in Feller's Eq. (27) by $\nu$ in order to achieve compatibility with the notation in \eqref{eq:BW-logistic}.

Feller starts out by calculating the  explicit solution to the KFE of  the pure linear death process,
that is, the number of inidividuals alive at time $t$; he states this
as the result of a recursive construction.
With a typo corrected  ($e^{\lambda t}-1$ must be replaced by $1-e^{-\lambda t}$ in his formula (6)), this same formula
says that
the number of inidividuals alive at time $t$
has a binomial distribution with parameters $N$ and $e^{-\lambda t}$
if there are initially $N$ individuals. (Today, after \cite{Feller 1950}, we would conclude this immediately, without solving systems of Kolmogorov forward equations,
via the probabilistic argument that there are initially  $N$ independent individuals,
each of
which dies at rate $\lambda$ and is therefore alive at time $t$
with probability $e^{-\lambda t}$.)
Likewise, the solution of the pure linear \emph{birth} process with
per capita birth rate $\lambda$, which he
gives in his Eq.~(17), is the negative binomial distribution
with parameters  $N$ and $e^{-\lambda t}$, which arises as the distribution of the sum
of $N$ independent random variables that are geometrically distributed
with parameter $e^{-\lambda t}$.
Again, this has a nice interpretation as the offspring of $N$ independently reproducing ancestors.

For the general birth process,
with arbitrary birth rates $p_n$, Feller notes that the KFE define
a probability distribution if and only if either only finitely many
of the $p_n$ are positive, or if $\sum_n 1/p_n$ diverges;
this is a standard textbook result today (usually presented in the
generalisation to birth-and-death processes). Under the conditions
stated, he also gives the explicit solution in passing.

\subsubsection{The moments of the stochastic process, and their
relationship with the deterministic equation}
Feller is particularly interested in the expectation, variance, and
other moments of the (random)
number of individuals alive at time $t$. In a trendsetting way, he
does not calculate them from the explicit solution of the KFE, even
where this is known; he rather uses the KFE to derive \emph{differential
equations for the moments}. Let $M(t)= \sum_k k P_k(t) = \mathbb E[K(t)]$ be the expected number of
individuals at time $t$. As stated above, Feller observes that, for the linear birth-and-death
process, $M(t)$ follows the differential equation for the deterministic
population model, and hence the expectation of the stochastic process
coincides with the deterministic solution.
In contrast, for the logistic
model, he finds from the differential equation relating the first and the second moment that
\begin{equation} \label{eq:BW-diffineq}
       \dot M < f(M)
\end{equation}
with $f$ of Eq.~\eqref{eq:BW-logistic}. From this he argues that $M$ is always less
than the  solution of the logistic equation. An alternative way to see
\eqref {eq:BW-diffineq} would be to observe that the KFE gives
$$\frac{d}{dt}\mathbb E[K(t)] = \mathbb E[g(K(t)]$$
where $g(k) = \sum_nQ(k,n)n$ \,\, and
$$Q(n,n+1)=\lambda n, \quad Q(n,n-1) = \gamma n^2, \quad Q(n,n) = -(\lambda n+ \gamma n^2), \quad Q(k,n) = 0 \text{\ \ otherwise.}$$
As a matter of fact, it turns out that $g(k) = \lambda k - \gamma k^2$, which is strictly concave, and hence \eqref{eq:BW-diffineq} is a consequence of Jensen's inequality.
Since Feller does not consider models with positive interaction (such as
symbiosis) in this part of the paper, he does not encounter the
convex situation.

\subsection{Diffusion equations for population dynamics}
\label{sec:BW-Depopdyn}
The second part of the paper (Sections 5--8 and 10) is devoted to the diffusion
limit of stochastic population dynamics. 
We cannot resist to quote Feller's thoughts from the beginning of Section 5, formulated in an almost literary German, about the \emph{substantially more lithesome probabilistic treatment, in which the population size is no longer assumed as integer-valued}, and where he alludes to similarities to the Brownian motion:
\begin{quote}\itshape
    Wir wenden uns nun der anderen von der in der Einleitung erw\"ahnten wahr\-scheinlichkeitstheoretischen Behandlungsweisen des Wachstumsproblems zu, welche we\-sent\-lich geschmeidiger ist, und bei der die Gr\"osse der Population nicht mehr ganzzahlig vorausgesetzt wird. Den Mechanismus des Vorgangs kann man sich hier \"ahnlich wie bei der Brownschen Molekularbewegung vorstellen. Der Zustand der betrachteten Population, d.\,h.\ ihre gesamte Lebensenergie ist einer dauernden Ver\"anderung unterworfen \textup{[\ldots]} 
\end{quote}

Starting from the transition density, Feller
calculates the infinitesimal drift $a(x)$ and  the infinitesimal variance
$b(x)$ (provided they exist).
With remarkable intuition, and a clear view of the branching property,
he states that, in the case of a stochastically independent reproduction
of the individuals, $a(x)$ and $b(x)$ must be proportional to $x$. Again, let us quote in German:
\begin{quote}\itshape
    Nimmt man beispielsweise an, dass die Gr\"osse der Population keinen Einfluss hat auf die durchschnittliche Vermehrungsgeschwindigkeit der Einzelindividuen, d.\,h.\ dass diese untereinander stochastisch unabh\"angig sind \textup{[\ldots]},  so m\"ussen $a(x)$ und $b(x)$ offenbar proportional zu $x$ sein \textup{[\ldots]}
\end{quote}
This gives rise to his equation (38), which reads as
\begin{equation}\label{eq:BW-Feller_branching}
\frac{\partial w(t,x)}{\partial t} =
\beta \frac{\partial^2 (x w(t,x))}{\partial x^2}
- \alpha \frac{\partial (x w(t,x))}{\partial x}.
\end{equation}
Here $w(t, \cdot)$ is the density of  population  size  at time
$t$, and $\alpha$ and $\beta$ are positive constants.
This seems to be the first appearance of what became famous as
\emph{Feller's branching diffusion}. We will come back to this in
Section~\ref{sec:BW-dpg}.

It is interesting to note that
the diffusion process in this second part of the paper is
\emph{not derived} from the birth-death jump processes
which Feller has presented in the first part; maybe, at this early stage,
the subtle rescaling required for this limit was not yet at his fingertips.
A decade later, however, he had these techniques; see
Section~\ref{sec:BW-dpg} on \emph{Diffusion processes in genetics}.
In 1939, Feller does allude to the birth-death processes,
but the connection is not yet clear. For example,
he tells us that the  $\alpha$ in \eqref{eq:BW-Feller_branching}
corresponds to the $\lambda$
encountered in the pure birth process. This is correct for the expected growth rate, but as a matter of fact
a pure birth process cannot have the diffusion limit
\eqref{eq:BW-Feller_branching}, since the paths of the former can only increase
in time, whereas the paths of the latter have fluctuations in both directions.

This issue reappears when Feller discusses the extinction probability of
the diffusion process. He  notes the important fact
that this quantity increases
with $\beta$, since it is tied to the fluctuations of the process, and at the same time emphasises as a sort of paradox
that, even for a positive net growth $\alpha > 0$, the diffusion process may die out with positive probability, while the population described by  the deterministic differential equation \eqref{eq:BW-linode}, as well as a pure birth process, cannot die out. (This paradox is resolved when one has in mind the different rescalings that lead to \eqref{eq:BW-linode} and \eqref{eq:BW-Feller_branching}.)

Following these considerations of the linear birth-and-death process,
Feller  includes dependence between individuals in Sections~8 and 10.
He presents two
specific examples. The first is his Equation (51), which is
the diffusion version of his Equation (7) and known today as
\emph{Feller's branching diffusion with logistic growth}
\cite{BW-Lambert2005,BW-LePardouxWakolbinger2013}. The second is the
(two-dimensional)
diffusion describing a two-species model with predator-prey
interactions, which now is also called \emph {Lotka-Volterra process}, see Eq. (1.2) in \cite{BW-CattiauxMeleard2010}. 
As with the jump processes in the first part of his paper, Feller is
concerned with the moments of the diffusion processes
and writes down a general recursion for
the $k$th moments $M_k$.
In the two-species model,
the interaction is positive (from the point of view of the
predator), so we finally encounter the convex case, in which
the expectation is greater than the solution of the corresponding
ODE.

In Section 9, Feller makes some
final remarks concerning the deterministic limit of both the
birth-death jump processes and the branching diffusion. These
are brief, heuristic calculations, which hint at the convergence
of the stochastic models to
Volterra's population models in the limit of infinite population size.
Today, powerful laws of large numbers
are available for large classes of such processes \cite[Chap.~11]{BW-EthierKurtz1986}.
They go far beyond the simple linear case alluded to in \eqref{eq:BW-linode}; rather, they
include quite general forms of density dependence. This leads us to the present
state of population dynamics.

\subsection{Afterthoughts}
Today, 75 years after \cite{Feller 1939a}, stochastic population dynamics
constitute a vibrant area of research, so wide that it is impossible
to give an overview in a short paragraph. Suffice it to say that
major questions raised by Feller continue to be ardent research themes.
Above all, this is true of interactions within and between populations.
Even simple models for the competition of two populations, whose
deterministic limit can be tackled as an easy exercise, turn into
hard problems when considered probabilistically.
Specifically, diffusion models with interaction
have become objects of intense research, see, e.g.,
\cite{BW-Etheridge2000,BW-Perkins2002} and references therein.

In the context of this commentary, it is particularly noteworthy
that a class of models known under the name of \emph{adaptive dynamics}
brings together ecological aspects (on a short time scale)
and genetical aspects (on a longer time scale) and thus builds
a bridge between population dynamics and population genetics.
A nice  overview of this topic and many others may be found
in the monograph by Haccou, Jagers, and Vatutin \cite{BW-HaccouJagersVatutin}.
Let us now turn to Feller's contributions to population
genetics.

\section{Feller and population genetics}
As already laid out in Section~\ref{sec:BW-intro}, important processes in population genetics are those that describe the evolution of type frequencies, or in other words, of proportions of subpopulation sizes  within a total population, whose size may vary as well. In this context, we may think of the individuals as \emph{genes}, where each gene is of a certain type, say $a$ or $A$.

The foundations of mathematical population genetics were laid starting in the 1920s by
Fisher, Wright, and Haldane. Their work mirrors the genetics of their time, today known as
\emph{classical genetics}. It had to rely on the phenotypic appearance of individuals (colour of
flower, surface structure of peas, body weight, milk yield \ldots). The molecular basis of
genetics was still unknown, so genes had to be treated as  abstract entities. When
\emph{molecular genetics} entered the labs in the 1960s, population genetics changed
dramatically, with Kimura as a leading figure, see Section~\ref{sec:BW-neutral}. The next (and, from
a 2014 perspective, the last) big leap took place in 1982, when Kingman introduced the
genealogical perspective via the \emph{coalescent process}. Comprehensive overviews of population genetics
theory are given in the textbooks by Ewens \cite{BW-Ewens2004} and
Durrett \cite{BW-Durrett2008}; for coalescent theory in particular, we further recommend
 Berestycki \cite{BW-Berestycki2009} (from a mathematical point of view) and
Wakeley \cite{BW-Wakeley2009} (from a more biological perspective).

With his contributions to population genetics, Feller thus was in the midst of
an important line of development.
We will comment on two of these articles.
The first, \emph{Diffusion processes in genetics} \cite{Feller 1951d},
is a landmark contribution towards stochastic modelling and
analysis via diffusion processes, and, as a matter of fact, reaches far
beyond population genetics as such.
The second, \emph{On fitness and the cost of natural selection} \cite{Feller 1967c}, uses
deterministic modelling (and is therefore similar in spirit to the
`Volterra equations').

\subsection{Diffusion processes in genetics}
\label{sec:BW-dpg}
Feller's article \emph{Diffusion processes in genetics} \cite{Feller 1951d} appeared in the Proceedings of the 2nd Berkeley Symposium on Mathematical Statistics and Probability, which took place in 1950. The central role of \cite{Feller 1951d} is nicely put into perspective by the following quote from Thomas Nagylaki's review \cite{BW-Nagylaki1989} on Gustave Mal\' ecot and the  transition from classical to modern population genetics:
\begin{quote}\itshape
    Mathematical research in diffusion theory influenced population genetics only gradually. As described in more detail below, Wright was unaware of Kolmo\-gorov's (1931) pioneering paper, and Wright, Mal\'ecot, and Kimura were all apparently unacquainted with Khintchine's (1933) book.\textup{[\ldots]}
    Thus, the mutually beneficial cross-fertilization between diffusion theory and population genetics did not start until Feller published his seminal 1951 paper.
\end{quote}

In the introduction of that paper, Feller sets the stage by writing:  
\begin{quote}\itshape
    Relatively small populations require discrete models, but for large populations it is possible to apply a continuous approximation, and this leads to processes of the diffusion type.
\end{quote}

Two diffusion processes are in the focus of the paper. One is what is  nowadays called \emph{Feller's branching diffusion}, the other is the so-called \emph{Wright--Fisher diffusion}.  Feller describes them by their \emph{diffusion equations}  (5.1) and (7.1), which are the Kolmogorov forward equations (or Fokker--Planck equations) for the densities, here called $u(t,x)$, cf. Section \ref{sec:BW-Depopdyn}.  Feller writes on pp.~228--229: 
\begin{quote}\itshape
    It is known that an essential part of Wright's theory is mathematically equivalent to assuming a certain diffusion equation for the \textup{gene frequency} \emph{(}that is, the proportion of $a$-genes\emph{)}.
\end{quote} 
In a footnote on the same page, Feller gives hints to the roots of this knowledge in the work of Kolmogorov, Fisher, Wright, and Mal\'ecot. 

\subsubsection{A foresight: Feller's diffusions as solutions of stochastic differential equations}\label{BW-twoSDEs}
Nowadays we do not hesitate to write the process described by Feller's equations (5.1) and (7.1) as solutions of stochastic differential equations in the sense of It\^o:
\begin{align}
\label{eq:BW-51}\tag{$\text{5.1}'$}  
    dZ_t &= \sqrt{2\beta Z_t} \, dW_t + \alpha Z_t \,dt,\\
\label{eq:BW-71}\tag{$\text{7.1}'$} 
    dY_t &= \sqrt{2\beta Y_t(1-Y_t)} \, dW_t + (\gamma_2 (1-Y_t)-\gamma_1Y_t) \,dt,
\end{align}
where $W$ is a standard Brownian motion. Feller legitimately resisted writing the processes in this form.  In \cite{Feller 1952c}, which grew out of Feller's invited lecture at the International Congress of Mathematicians in the year 1950, he writes about It\^o's Stochastic Analysis:  \begin{quote}\itshape
    This approach has the advantage that it permits a direct study of the properties of the path functions, such as their continuity, etc. In principle, we have here a possibility of proving the existence theorems for the partial differential equations \textup{[\ldots]} directly from the properties of the path functions. However, the method is for the time being restricted to the infinite interval and the conditions on \textup{[the diffusion and drift  coefficients]} $a$  and $b$ are such as to guarantee the uniqueness of the solution. So far, therefore, we cannot obtain any new information concerning the ``pathological'' cases.
\end{quote}

\subsubsection{An emerging theme: What happens at the boundary?}
Indeed, the state spaces of (5.1) and (7.1) are not the `infinite interval' $(-\infty, \infty)$ but  $[0, \infty)$ and $[0,1]$,  and it took 20 years until T.\ Yamada and S.\ Watanabe proved that the coefficients in the above stated SDEs are good enough to guarantee strong uniqueness of the solution,  see \cite{BW-YamadaWatanabe1971} and also \cite{BW-Watanabe1969}. A coupling argument from Stochastic Analysis then guarantees that the solution of \eqref{eq:BW-71} converges in law to the unique equilibrium distribution whose density is the unique invariant probability density of (7.1), which is the  Beta($\gamma_2/\beta, \gamma_1/\beta$)-density. Thus, although  for $\gamma_1 < \beta$  the random path $Y$ hits $0$ with probability one (and similarly for $\gamma_2 < \beta$ it hits $1$ with probability one), these visits to the boundary do not lead (as conjectured by Feller on p.~239) to a non-vanishing \emph{accumulation of the masses concentrated at $x=0$ and $x=1$ which is maintained in the steady state}, in other words, the coefficient $\mu$ in his equation (7.3) is in fact equal to 1.

Questions like these may have been one source of motivation for Feller to initiate his groundbreaking studies on the boundary classification of diffusion processes, see his footnote on p.~234, where he speaks of \emph{boundary conditions of an altogether new type}, and the one on p.~229 \emph{added in proof}, where he announces that \emph{a systematic theory, including the new boundary condition, is to appear in the Annals of Math.} Feller's \emph{classification of boundaries}  is reviewed and commented in Section 2 of Masatoshi Fukushima's essay in this volume.

\subsubsection{The diffusion approximation of the Wright--Fisher chain and beyond}
As already indicated, another important aspect that is taken up in Feller's paper is that of the \emph{diffusion approximation}, i.\,e.~the convergence of a sequence of (properly scaled)  discrete processes to the solution of (5.1) and (7.1), respectively. In the former case the underlying discrete process is a Galton--Watson process, in the latter it is the Wright--Fisher Markov chain. The transition probabilities of the Wright--Fisher  chain are given by (3.2), (3.4) and (3.5).  The  \emph{diffusive mass-time-scaling} is given by (8.5):  a unit of time consists of $N$ generations,  and a unit of mass consists of $N$  (or here $2N$) individuals, with $N$ being the total population size.   The scaling (8.4) of the individual mutation probabilites $\alpha_1, \alpha_2$ is that of \emph{weak mutation}, which leads in the scaling limit to the \emph{ infinitesimal mean displacement} $a(x)$ and the \emph{infinitesimal variance} $2b(x)$, see  (8.6) and (8.7).  In the context of (7.1), the \emph{ drift coefficient} $a(x)$ is due to the effect of mutation, and the \emph{diffusion coefficient} $b(x)$ describes the strength of the fluctuations that come from the random reproduction.  (In order to be consistent with (7.1), $\beta_i$ should be replaced by $\gamma_i$ in  (8.4), (8.6) and (8.9)). The `convergence of generators' which emerges from (8.6) and (8.7)  can be lifted to the convergence of the corresponding semigroups, see e.\,g.\ the chapter on \emph{Genetic Models} in the  monograph \cite{BW-EthierKurtz1986} by Ethier and Kurtz.

The convergence theorems in \cite{BW-EthierKurtz1986} comply with Feller's programmatic proposal: \emph{It should be proved that our passage to the limit actually leads from \textup{(8.2)} to \textup{(8.6)}}, i.\,e.\ from the probability weights of the Wright--Fisher chain to the probability densities the Wright--Fisher diffusion. To achieve this, Feller proposed an  expansion into eigenfunctions (in particular he found the eigenvalues of the Wright--Fisher transition semigroup) and checked part of the convergence in Section 8 and Appendix I.  Such a representation is not required in the systematic approach presented in \cite{BW-EthierKurtz1986}.
Still, the approach via eigenfunctions is interesting in its own right, and has been extensively used in Mathematical Biology.

At the beginning of  Section 9 (entitled \emph{Other possibilities}) Feller writes :
\begin{quote}\itshape 
    The described passage to the limit which leads to Wright's diffusion equation \textup{(7.1)} is different from the familiar similar processes in physical diffusion theory where the ratio $\Delta x/\Delta t$ tends to infinity rather than to a constant. It rests entirely on the assumption \textup{(8.4) [of weak mutation]}. We shall now see that \textup{any modification of this assumption leads to a non-singular diffusion equation of the familiar type \textup{(}to normal distributions\textup{)}}.
\end{quote}
Indeed, for the scaling (9.1), (9.2)  $\alpha_i = \gamma_i\epsilon,\, i=1,2$, $N\varepsilon \to \infty$, which corresponds to \emph{strong mutation}, Feller states a law of large numbers, i.\,e.~a convergence of the type frequencies to the equilibrium point $\left(\frac{\gamma_2}{\gamma_1+\gamma_2}, \frac{\gamma_1}{\gamma_1+\gamma_2}\right)$, and argues that the (properly scaled) process of fluctuations around this equlibrium point converges to a  process whose probability density satisfies the diffusion equation (9.10) (and thus is an Ornstein-Uhlenbeck process).
\subsubsection{The diffusion approximation of Galton--Watson processes}
The diffusion equation (5.1) appeared already in \cite{Feller 1939a}, see Eq.~\eqref{eq:BW-Feller_branching} in Section \ref{sec:BW-Depopdyn}. However, as we have
seen there, certain issues concerning the (scaling) limits of Galton--Watson processes had remained unrevealed in  \cite{Feller 1939a}. Towards 1950, Feller was ready to  attack this. As to the diffusion approximation of a sequence of ``nearly critical'' Galton--Watson processes by (5.1), Feller gives a proof in Appendix II. His idea is to take the iterates $f_n$  of the offspring generating function (which are known to describe the generating functions of the subsequent generation sizes) to their scaling limit. This limit turns out to satisfy the PDE (12.9) (which, in turn, corresponds to (5.1)). Feller writes: 
\begin{quote}\itshape
    We effect this passage to the limit formally: it is not difficult to justify these steps, since the necessary regularity properties of the generating functions $f_n(x)$ were establihed by Harris \textup{\cite{BW-Harris1948}}.
\end{quote}  
Again, from today's perspective, an alternative way is provided by the convergence of generators, see \cite{BW-EthierKurtz1986}. In the very last lines of his Appendix II, Feller remains a bit sketchy when he writes that 
\begin{quote}\itshape
    the boundary condition $u(t,0)$ follows from the fact that in the branching process the probability mass flowing out into the origin tends to zero.  
\end{quote}
In fact, for the solution $Z$ of \eqref{eq:BW-51}  (with $Z_0 = 1$, say),  the probability mass flowing out into the origin is non-zero at any fixed time $t $, and the density of $Z_t$ does not vanish near the origin. Again, the desire to obtain clarity  on questions like these may have been a motivation for Feller's then upcoming research on the boundary behaviour of one-dimensional diffusions.

\subsubsection{From two-type to continuum-type generalizations of Feller's branching diffusion}
In the introduction, Feller points out that \emph{serious difficulties arise if one wishes to construct population models with interactions among the individuals}, and that \emph{the situation grows worse if the population consists of different types of individuals}. He then writes: \begin{quote}\itshape 
    In fact, the bivariate branching process leads to such difficulties that apparently not one single truly bivariate case has been treated in the literature. In the theory of evolution this difficulty is overcome by the assumption of a constant population size \textup{[\ldots]} In Section 10 the assumptions of constant population size is dropped and a truly bivariate model is constructed which takes into account selective advantages in a more flexible way. \textup{[\ldots]} The same limiting process which leads \textup{[\ldots]} to the diffusion equation of Wright's theory can be applied to our new bivariate model and leads to a diffusion equation in two dimensions.
\end{quote} 
These \emph {two-dimensional Markov processes with branching property} have been taken up and analysed in a broader context in 1969 in the paper \cite{BW-Watanabe1969} which carries that title.  Already before, Watanabe had published his seminal paper \cite{BW-Watanabe1968} which established Feller's branching diffusions with a continuum of types.  This together with the poineering work of Don Dawson 
gave rise to a class of processes that were later called Dawson--Watanabe superprocesses (\cite{BW-Etheridge2000, BW-DawsonPerkins2012}). A good part of Perkins' Saint Flour Lecure Notes (part 2 of \cite{BW-DawsonPerkins2012}) is devoted to superprocesses with interactions, and thus is fully on the line of Feller's program to construct population models with different types of individuals and with interactions among the individuals.

\subsubsection{The inner life of Feller's branching diffusion: excursions and continuum trees}
This is a good place to mention another fascinating development which is connected with Feller's branching diffusions and is associated with the names of  Daniel Ray and Frank Knight (the latter was Feller's doctoral student and Ed Perkins' PhD advisor).

Thanks to the branching property (and the thereby implied infinite divisibility), the random path $Z$ of a Feller branching diffusion is a Poissonian sum of countably many ``Feller excursions'' $\zeta$. In fact, each of them has an ``internal life'' in the sense that $\zeta_t $ is the size at time $t$ of a continuum population originating from one single ancestor. The genealogical tree of this population can be described by a Brownian excursion $\eta$ reaching level $t$, which can be imagined as the `exploration path' of a continuum random tree whose mass alive at level $t$ is $\zeta_t$.   The second  Ray--Knight theorem says that $\zeta_t $ can be represented as the local time spent by $\eta$ at level $t$. In a discrete setting of Galton--Watson processes, this is ancticipated in Harris' work \cite{BW-Harris1952} with its section on \emph{walks and trees}. The correspondence between a Feller excursion $\zeta$ and an It\^o excursion $\eta$    is depicted on the first page of \cite{BW-PardouxWakolbinger2011}, framed by pictures of  Feller and It\^o, who met in person at Princeton in 1954. See \cite{BW-PardouxWakolbinger2011} for more explanations, and references to groundbreaking developments that dealt with the genealogical  structure behind Feller's branching diffusion, such as Aldous' \emph{Continuum Random Tree} (which plays in the world of random trees a similar role to that of Brownian motion in the classical invariance principle) and  Le Gall's \emph{Random Snake}, which  provides a representation of the Dawson--Watanabe super-Brownian motion as a continuum-tree-indexed Markov motion.

For more on excursions and excursion point processes in relation with Feller's work, see Section 3.1 of the contribution of M.~Fukushima.
\subsubsection{Frequencies in multivariate continuum branching: conditioning and time change}
Another interesting question which Feller addresses at the end of his introduction concerns the relative frequencies in a bivariate model of branching diffusions. Feller writes: 
\begin{quote}\itshape \textup{[\ldots]} 
    it is to be observed that \textup{in no truly bivariate case does the gene frequency satisfy a diffusion equation} \emph{(}Sections 6 and 10\emph{)}. In fact, if the population size is not constant, then the gene frequency is not a random variable of a Markov process. Thus, conceptually at least, the assumption of a constant population size plays a larger role than would appear on the surface.
\end{quote}

Indeed, as it turns out (and Feller may have been well aware of this), one way of passing from (5.1) to (7.1), say with $\alpha = \gamma_1=\gamma_2=0$, is to consider two independent Feller branching diffusions (solutions of \eqref{eq:BW-51}) $Z^{(1)}$ and $Z^{(2)}$, conditioned to $Z^{(1)}+ Z^{(2)}= 1$. Of course, this must be given a precise meaning, and this has been done in a much more general context in the papers by Etheridge and March \cite{BW-EtheridgeMarch1991} and Perkins \cite{BW-Perkins1991}. The title of Perkins' paper is programmatic: A Dawson--Watanabe superprocess conditioned to have constant mass one is a Fleming--Viot process (which is the continuum-type, and thus measure-valued, generalization of the Wright--Fisher diffusion). In the context of (5.1) to (7.1) this means that under the conditioning $Z^{(1)}+ Z^{(2)}= 1$ the process  $Y:= Z^{(1)}$ is a Wright--Fisher diffusion. On the level of genealogies, the conditioning to a constant total mass takes the continuum random forest that underlies \eqref{eq:BW-51} into Kingman's coalescent.

A second way to get from \eqref{eq:BW-51} to \eqref{eq:BW-71} (again for $\alpha = \gamma_1=\gamma_2=0$) is to consider the relative frequency $Y= \frac{Z^{(1)}}{Z^{(1)}+ Z^{(2)}}$ after a time change $ds = dt/(Z^{(1)}+ Z^{(2)})$. In this way, the relative frequencies again become Markovian and, by an easy application of It\^o's formula, turn out to solve \eqref{eq:BW-71}.

\bigskip
The concept of time change is central also in the work of John Lamperti. Lamperti's work can be seen as a direct continuation and extension of Feller's ideas, introducing and analysing the continuum mass limits of Galton--Watson processes also for heavy-tailed offspring distributions \cite{BW-Lamperti1967a}. His  article \cite{BW-Lamperti1967b}, which was communicated by H.\,P.~McKean, another former PhD student of Feller, introduced what is now called Lamperti's transform, a time change which establishes the link between L\'evy processes and continuous state branching processes.

\bigskip
To conclude: Feller's paper \emph{Diffusion processes in genetics}  is a remarkable contribution at the interface of probability theory and population biology, with enduring stimulations in either direction. It takes a central place  in the development of  mathematical population genetics, and has triggered substantial new directions in the modern theory of stochastic processes.

\subsection{The cost of natural selection}
Feller's article \emph{On fitness and the cost of natural selection} \cite{Feller 1967c} appeared
in \emph{Genetical Research Cambridge}, a renowned biology journal. The
introduction contains the disclaimer 
\begin{quote}\itshape
This paper is written by
a mathematician, and accordingly no new biological models or hypotheses are
advanced.
\end{quote} 
It may be added that the mathematics is fairly elementary
as well, and the strength of the article is the concise conceptual
thinking by which Feller dissects the logics of an argument
that enormously influenced the genetic thinking of that time,
and finds a fundamental weak point in it.

\subsubsection{The 1960s and the neutral theory of population genetics}
\label{sec:BW-neutral}
As hinted at already, the 1960s were turbulent times for genetics---and for population
genetics in particular. Until then,
 the variation between individuals
could only be observed at the phenotypical level, and much of this
was easily explained by selection: Stronger beaks crack harder nuts;
webbing eases swimming; fat pads protect against the cold. Then,
in the 1960s, the first observations of \emph{variation at the molecular
level} became available -- not yet via sequencing, but  via
so-called restriction length polymorphsims (RFLP) of DNA or via
gel electrophoresis of proteins. The resolution of these methods
is lower than that of sequencing but, nevertheless, the variation was
so much larger than expected on phenotypic grounds that researchers
were shocked---and were puzzled about the question: \emph{Can this all
be explained by selection?}

These considerations were strongly influenced by the concept
of the \emph{genetic load}, coined by Haldane \cite{BW-Haldane1957};
in particular the concept of the \emph{substitutional
load}. This is the number of \emph{selective deaths}, that is,
the number of individuals `killed' by selection in the process
of substituting  one type by a fitter one. Put differently, the
substitutional load (or \emph{cost of natural selection})
is the number of excess individuals that must be produced in a population
under selection. In 1968, Kimura \cite{BW-Kimura1968} concluded that, if a large
fraction of the observed variability is selective,
the load is forbidding. This led to one of the most influential
and conflict-prone hypotheses of evolutionary theory, namely, to the
so-called \emph{neutral theory}, which claims that the overwhelming
proportion of the observed molecular variation is selectively neutral,
that is, most mutations do not change fitness.

In what follows, we look more closely into the concept of the substitutional
load, and into Feller's criticism of it. We restrict ourselves to the
case of haploid populations (i.e., carrying only one copy of
the genetic information) that reproduce asexually (Sections 1--5 in
Feller's paper). In Sections 6--9, he tackles additional complications
that emerge in diploid individuals (with two copies of every gene),
but the conceptual issues  are more transparent in the haploid case.

\subsubsection{Absolute and relative frequencies in population genetics}
Feller considers a population of individuals that consists of our two
types $A$ and $a$, large enough to justify deterministic treatment.
He assumes discrete generations  where every $A$-individual 
leaves an average of $\mu$ offspring for the next generation, whereas every $a$-individual
produces  an average of $\mu'=\mu (1-k)$ descendants, $0 < k < 1$.
The quantities $\mu$ and $\mu'$ are known as (absolute) fitnesses in population genetics.
Each of the two
subpopulations then grows (or decays) geometrically,
\begin{equation}\label{eq:BW-absfrequ}
   N_n = \mu N_{n-1}, \quad N_n' = \mu (1-k) N_{n-1}',
\end{equation}
so $N_n= N_0 \mu^n$ and
$N_n'= N_0 '(\mu (1-k) )^n$, where $N_n$ and $N_n'$ denote the size of
the $A$- and $a$-subpopulations, respectively.

Now population genetics is traditionally more concerned with
relative frequencies of types than with absolute ones; the main
reason is that relative frequencies are easier to measure. One
therefore considers
\[
  p_n := \frac{N_n}{N_n+N_n'}, \quad q_n := \frac{N_n'}{N_n+N_n'}.
\]
Clearly, under \eqref{eq:BW-absfrequ},
\begin{equation}\label{eq:BW-relfrequ}
  q_n = \frac{q_0 (1-k)^n}{p_0 + q_0 (1-k)^n},
\end{equation}
which is Feller's Eq.~(3.5).
Obviously, the powers of $\mu$ cancel out, which is a strength and a
weakness at the same time: On the one hand, this means that knowledge
only of the ratio of the fitness values is required to predict
the behaviour of the population. Actually, \eqref{eq:BW-relfrequ}
holds more generally than the simple derivation may suggest: It continues to
hold if \eqref{eq:BW-absfrequ} is replaced by
\begin{equation}\label{eq:BW-eco}
  N_n = \mu N_{n-1} f(N_n+N_n'), \quad N_n' = \mu (1-k) N_{n-1}' f(N_n+N_n').
\end{equation}
Here $f$ is a function that depends on the total population size only
and acts on both types in the same way.
In typical ecological models, one uses a  monotonically decreasing function $f$
with $f(0)=1$ in order to describe
how population size decreases the
per capita offspring size
due to competition. In particular, for suitable choices of $f$,
both the population size and the relative frequencies will, in
the long run, approach stationary values.

The downside of thinking in terms of relative fitnesses is that
one loses sight
of the absolute population sizes. The latter may lead to absurd
conclusions, in particular in cases where one or both subpopulations
go extinct. This leads us to Feller's criticism of the genetic load.

\subsubsection{The substitutional load and Feller's criticism}
Feller now recalls Haldane's definition of the genetic load:
In generation $n$, the mixed population has a loss of $d_n := k q_n$
offspring relative to a population consisting of $A$ individuals only,
where $d_n$ must be measured in units of the total population size
of generation $n$.  Over  $M$ generations, Haldane takes
$D:= \sum_{n=0}^Md_n$ as the total genetic load. Calculating this
with reasonable parameters and as $M \to \infty$, he arrives at a representative value
of 30 for the cost of substitution of one gene by a fitter one.

Haldane's definition of $D$ makes sense
only if population size remains constant over time, or if changes
in population size are so small that they can safely be neglected.
Haldane does not make this explicit; Feller has an eye on
both possibilities.

Feller first considers the case that the population size is not
constant; rather, $N_n$ and $N_n'$
behave as  in \eqref{eq:BW-absfrequ}. Since $\mu=1$ and $k<1$ is assumed,
this means that $N_n \equiv N_0$ and $N_n' \rightarrow 0$ as $n \to
\infty$, so the total population size decreases from $N_n+N_n'$ to
$N_0$. On an absolute scale, the loss of individuals due to selection is
$N_n'-N_{n+1}' = k N_n'$ in generation $n$, and altogether
\[
(N_0'-N_1') + (N_1'-N_2') + \ldots = N_0',
\]
in agreement with the
decrease of the total population size from $N_0+N_0'$ to $N_0$.
In contrast, Haldane's $D$, which neglects the size change,
can give much larger values; in particular,
it can be larger than the number of $a$ individuals ever born. Obviously,
this understanding of the load produces severe artifacts.

Feller then discusses how Haldane's argument may be `rescued'.
One possibility is to keep population size constant by immigration
from a reservoir population, in exactly the same proportions
as in the current population under consideration. Then, Haldane's
$D$ gives the correct answer
(but it must be kept in mind that the cost of selection is then borne
by the reservoir population). The second possibility is to
assume an ecological model rather than geometric growth. Feller speaks
of $\mu$ depending on population size; maybe as in our
Eq.~\eqref{eq:BW-eco}, where $\mu$ is replaced by $\mu f(N_n+N_n')$.
But more general models are also possible; for example, each type
may be affected by competition in a different way. The genetic
load would then be identified with the decrease in the stationary
population size due to the reduced reproduction rate of one type.
But now there is a wealth of possible models, and the genetic
load would depend on the details. Specifying these is a
task Feller assigns to the biologists.

\subsubsection{Feller's criticism in the general context}
With remarkable insight, Feller dissects a conceptual problem of his time:
Load arguments can be inconsistent if they blindly rely on relative frequencies.
His fellow researchers in biology do, however, not seem to have taken
too much notice of his criticism. After all, as already mentioned
in Section~\ref{sec:BW-neutral}, a year after Feller's paper (and without
citing it), Kimura did
put forward his neutral theory of molecular evolution, to a large
extent on the basis of load arguments \cite{BW-Kimura1968}.
More precisely, Kimura
used an extension of Haldane's argument. He
assumed a sequence of numerous genetic loci (rather than Haldane's
and Feller's  single locus),  each of which
can be of a favourable or a less favourable type, with fitness
assumed as \emph{multiplicative across loci}.
As a consequence, there is a multitude of possible genotypes;
fitness differences between individuals can become enormous;
and the load
(if calculated in Haldane's manner) can become astronomical
(Kimura and Ohta~\cite{BW-KimuraOhta1971} give a  value of $D=10^{78}$).
Here the misconception lies
in the assumption of multiplicativity over loci, which is completely
unrealistic, but was hardly questioned at that time.
For the details, see the insightful presentation in
\cite[Chapter~2.11]{BW-Ewens2004}.

Let us return to the original question that load arguments
were supposed to answer: Can all the variation observed at the
molecular level be explained by selection?
Indeed, today, a large fraction of the molecular variation is
considered selectively neutral or nearly so, although single
mutations with spectacular selective effects are  well known.
But this insight is no longer built on load arguments  --
rather, the assumption of neutrality has proved extremely
successful in describing patterns of genetic variation.

Last but not least, it should be noted that there is a lot of
truth in Feller's general warning
not to neglect population size in population genetics. Indeed,
load arguments are not the only artifacts of this kind. Another example is
the famous phenomenon of \emph{Muller's ratchet}, which describes the
\emph{ad infinitum} accumulation of deleterious mutations due to
stochastic effects in finite populations of constant size.
If described in terms of an ecologically more realistic
(and conceptually more correct)
model with variable population size,  the accumulation does not
continue forever.
Rather, when fitness has declined below a threshold value, the
population  experiences a \emph{mutational meltdown}, which
ultimately leads to extinction (see \cite{BW-Baake-Gabriel2000} for a review).

\subsubsection*{Acknowledgement}
The interaction between mathematics and biology, which was greatly
advanced by Feller, has a continuing impact on both fields.
The authors gratefully acknowledge support from the Priority
Programme \emph{Probabilistic Structures in Evolution \emph{(}SPP 1590\emph{)}}
funded by Deutsche Forschungsgemeinschaft (German Research Foundation, DFG).

\end{document}